\documentclass[12pt,twoside,reqno]{amsart}
\usepackage[colorlinks=true,citecolor=blue]{hyperref}
\usepackage{mathptmx, amsmath, amssymb, amsfonts, amsthm, mathptmx, enumerate, color,enumitem}
\usepackage{graphicx}
\usepackage{multirow}
\usepackage{epstopdf}
\usepackage{multicol}
\usepackage{algorithm}
\usepackage{algorithmic}
\usepackage{booktabs}

\usepackage{bigstrut}
\usepackage[figuresright]{rotating}
\usepackage{adjustbox}
\graphicspath{{figures/}}
\usepackage{subfigure}
\usepackage[american]{babel}
\usepackage{microtype}
\usepackage{float}
\numberwithin{algorithm}{section}

\newtheorem{theorem}{Theorem}[section]

\newtheorem{lemma}{Lemma}[section]

\theoremstyle{definition}

\newtheorem{example}{Example}[section]
\newtheorem{remark}{Remark}[section]

\numberwithin{equation}{section}
\begin{document}
\setcounter{page}{1}

\vspace*{0.9cm}
\title[two inertial extragradient-mann algorithms]
{Self adaptive inertial extragradient algorithms for solving variational inequality problems}
\author[B. Tan, J. Fan, S. Li]{Bing Tan, Jingjing Fan,  Songxiao Li$^{*}$}

\maketitle
\vspace*{-0.6cm}

\begin{center}
{\footnotesize {\it
 Institute of Fundamental and Frontier Sciences,\\ University of Electronic Science and Technology of China, Chengdu,  China\\
 
}}
\end{center}

\vskip 4mm {\small \noindent {\bf Abstract.}
In this paper, we study the strong convergence of two Mann-type inertial extragradient algorithms, which are devised  with a new step size, for solving a variational inequality problem  with a monotone and Lipschitz continuous operator in real Hilbert spaces.
Strong convergence theorems for our algorithms are proved without the prior knowledge of the Lipschitz constant of the operator. Finally, we provide some numerical experiments to illustrate the performances of the proposed algorithms and provide a comparison with related ones.

\noindent {\bf Keywords.}
Variational inequality problem; Subgradient extragradient algorithm; Tseng's extragradient algorithm; Inertial algorithm; Mann-type algorithm.}

\renewcommand{\thefootnote}{}
\footnotetext{ $^*$Corresponding author.
\par
\hspace{-10pt}Email addresses: bingtan72@gmail.com (B. Tan),  fanjingjing0324@163.com (J. Fan), jyulsx@163.com (S. Li).}

\section{Introduction}
Let $C$ be a convex and closed set in a  real Hilbert spaces $H$ with the inner product $\langle\cdot, \cdot\rangle$ and the   norm $\|\cdot\|$. For all $ x,y\in H $, one recalls that a mapping $ T: H \rightarrow H $ is said to be (i) $L$-Lipschitz continuous with $L>0$ iff $\|T x-T y\| \leq L\|x-y\|$ (if $ L=1 $, then $ T $ is said to be  nonexpansive); (ii) $\eta$-strongly monotone if there exists $\eta>0$ such that $ \langle T x-T y, x-y\rangle \geq \eta\|x-y\| $; (iii) monotone if  $\langle T x-T y, x-y\rangle \geq 0$. A point $x^{*} \in H$ is called a fixed point of $T$ if $T x^{*}=x^{*}$. The set of all the fixed points of $T$ is denoted by $\operatorname{Fix}(T)$.  Let $A: H \rightarrow H$ be an operator. The variational inequality problem (shortly, VIP) for $A$ on $C$ is to find a point $x^{*} \in C$ such that
\begin{equation*}\label{VIP}
\left\langle A x^{*}, x-x^{*}\right\rangle \geq 0,\quad \forall x \in C\,.
\tag{\text{VIP}}
\end{equation*}
From now on, the solution set of \eqref{VIP} is denoted by $\mathrm{VI}(C, A)$.

In a wide range of applied mathematical problems, the existence of a solution is equivalent to the existence of a solution to the above classical variational inequality.  Therefore, the variational inequality   is an important tool in studying a wide class of physics, engineering, economics and optimization theory, see, e.g., \cite{IY1,IY2,QA,ANQ}. Over the last 60 years or so, the   variational inequality has been revealed as a very powerful and important tool in the study of various linear and nonlinear phenomena. Some problems, such as, systems of equations, complementarity problems, and   equilibrium problems, can be formulated as the variational inequality.

Recently, many authors   proposed and investigated various algorithms for solving the variational inequality, see, e.g., \cite{PRGM,fanjnca,dongopt,liujnca,fanopt,iPC} and the references therein. Projection methods and their variant forms act as important tools for finding  approximate solutions of the variational inequality. One of  well-known solution methods for \eqref{VIP} is the
following projection gradient algorithm:
\begin{equation}\label{PC}
\begin{aligned}
x_{n+1}=P_{C}\left(x_{n}-\lambda A x_{n}\right), \quad \forall n \geq 1\,,
\end{aligned}
\end{equation}
where $\lambda$ is a positive real number and $P_C$ is the metric (nearest point) projection onto $C$. However, the convergence of the algorithm requires a strongly monotonicity on $A$ (or inverse strongly, which is also usually said to be  cocoercive).
If $A$ is $L$-Lipschitz continuous and monotone, Korpelevich~\cite{EGM} proposed the following extragradient algorithm with double projections to reduce the monotonicity of operator $A$:
\begin{equation}\label{EGM}
\left\{\begin{aligned}
&y_{n}=P_{C}\left(x_{n}-\lambda A x_{n}\right)\,, \\
&x_{n+1}=P_{C}\left(x_{n}-\lambda A y_{n}\right), \quad \forall n\geq 1\,,
\end{aligned}\right.
\end{equation}
where $\lambda\in (0,\frac{1}{L})$. The algorithm converges to an element of $\mathrm{VI}(C, A)$ provided that $\mathrm{VI}(C, A)$ is non-empty.
In fact, in \eqref{EGM}, the price is that one needs to calculate two projections from $H$ onto the feasibility set $C$. If $C$ is a general convex-closed set, this might require a prohibitive amount of computation time. To overcome this computational drawback, many authors have modified this method in various ways. Next, we introduce two modifications of the extragradient algorithm.

The extragradient algorithm was modified by Tseng~\cite{tseng} with a remarkable scheme. The Tseng's extragradient algorithm reads as follows:
\begin{equation}\label{Tseng}
\left\{\begin{aligned}
&y_{n}=P_{C}\left(x_{n}-\lambda A x_{n}\right), \\
&x_{n+1}=y_{n}-\lambda\left(A y_{n}-A x_{n}\right), \quad \forall n\geq 1,
\end{aligned}\right.
\end{equation}
where $\lambda\in (0,\frac{1}{L})$. In 2011, Censor et al.~\cite{SEGM} modified the extragradient algorithm by replacing the second projection onto the  convex and closed subset with the one onto a subgradient half-space. The subgradient extragradient algorithm is of the form:
\begin{equation}\label{SEGM}
\left\{\begin{aligned}
&y_{n}=P_{C}\left(x_{n}-\lambda A x_{n}\right)\,, \\
&T_{n}=\left\{x \in H \mid \left\langle x_{n}-\lambda A x_{n}-y_{n}, x-y_{n}\right\rangle \leq 0\right\}\,, \\
&x_{n+1}=P_{T_{n}}\left(x_{n}-\lambda A y_{n}\right), \quad \forall n\geq 1\,,
\end{aligned}\right.
\end{equation}
where $\lambda\in (0,\frac{1}{L})$. We point out here that the Tseng's extragradient algorithm and the subgradient extragradient algorithm only need to calculate one projection onto $C$ in each iteration. Note that under some appropriate settings, Algorithm~\eqref{Tseng} and Algorithm~\eqref{SEGM}  converge to the solution of the variational inequality weakly. For this reason, a natural question that arises is how to design an algorithm, which provides strong convergence to solve problem \eqref{VIP} , when $A$ is only $L$-Lipschitz continuous and monotone mapping. Recently, Kraikaew and Saejung~\cite{HSEGM} based on  the subgradient extragradient algorithm and the Halpern algorithm  to proposed an algorithm for solving \eqref{VIP}. Their algorithm is of the form:
\begin{equation*}\label{HSEGM}
\left\{\begin{aligned}
&y_{n}=P_{C}\left(x_{n}-\lambda Ax_{n}\right)\,, \\
&T_{n}=\left\{x \in {H}\mid\left\langle x_{n}-\lambda Ax_{n}-y_{n}, x-y_{n}\right\rangle \leq 0\right\}\,, \\
&x_{n+1}=\alpha_{n} x_{0}+\left(1-\alpha_{n}\right) P_{T_{n}}\left(x_{n}-\lambda Ay_{n}\right),\quad \forall n\geq 1\,,
\end{aligned}\right.
\tag{\text{HSEGM}}
\end{equation*}
where $\lambda\in (0,\frac{1}{L}),$ and $ \alpha_n \subset (0,1)$ with  $\sum^{\infty}_{n=1}\alpha_n=+\infty$ and $\lim_{n\rightarrow\infty}\alpha_n=0$. They proved that $\{x_n\}$ converges   to the unique solution of \eqref{VIP} in norm. Note that the algorithm \eqref{HSEGM} needs to know the Lipschitz constant of  $A$, which limits the applicability of the algorithm. To overcome this shortcoming, Yekini and Olaniyi \cite{shefu} proposed a modification of the subgradient extragradient algorithm with the adoption of the Armijo-like step size rule. Indeed, they investigated the following algorithm:
\begin{equation*}\label{VSEGM}
\left\{\begin{aligned}
&\text{Given }\ell \in(0,1),\,\mu \in(0,1)\,,\\
&y_{n}=P_{C}\left(x_{n}-\lambda_{n} A x_{n}\right)\,, \\
&T_{n}=\left\{x \in H\mid\left\langle x_{n}-\lambda_{n} A x_{n}-y_{n}, x-y_{n}\right\rangle \leq 0\right\}\,, \\
&z_{n}=P_{T_{n}}\left(x_{n}-\lambda_{n} A y_{n}\right)\,,\\
&x_{n+1} = \alpha_{n}f(x_{n})+(1-\alpha_{n})z_{n},\quad \forall n\geq 1\,,
\end{aligned}\right.
\tag{\text{VSEGM}}
\end{equation*}
where $f: H \rightarrow H$ is a contraction mapping, $\lambda_{n}=\ell^{m_{n}}$ and $ m_{n} $  is the smallest nonnegative inter such that $\lambda_{n}\left\|Ax_{n}-Ay_{n}\right\| \leq \mu\left\|x_{n}-y_{n}\right\|$. They proved that the algorithm defined by \eqref{VSEGM} converges to the solution set of \eqref{VIP} in norm. This algorithm is not required to know the Lipschitz constant of   $ A $, but the step size needs to calculate the value of  $ A $ many times at each iteration. Therefore, although the Armijo criterion may not need to know the Lipschitz constant, it is very computationally expensive. Recently, Yang and Liu~\cite{TVEGM} combined the Tseng's extragradient algorithm and the viscosity algorithm with a simple step size and proposed a new iterative algorithm. The algorithm  consists of  only one projection and does not require the prior knowledge of the Lipschitz constant of the operator. They obtained a strong convergence theorem under suitable conditions, and their algorithm is described as follows:
\begin{equation*}\label{TVEGM}
\left\{\begin{aligned}
&\text{Given }\lambda_{0} \in(0,1),\, \mu \in(0,1)\,,\\
&y_{n}=P_{C}\left(x_{n}-\lambda_{n} A x_{n}\right)\,, \\
&z_{n}=y_{n}-\lambda_{n}(Ay_{n}-Ax_{n})\,,\\
&x_{n+1} = \alpha_{n}f(x_{n})+(1-\alpha_{n})z_{n}\,,\\
&\lambda_{n+1}=\left\{\begin{array}{ll}
\min \left\{\frac{\mu\left\|x_{n}-y_{n}\right\|}{\left\|Ax_{n}-Ay_{n}\right\|}, \lambda_{n}\right\}, & \text { if } Ax_{n}-Ay_{n} \neq 0\,; \\
\lambda_{n}, & \text { otherwise}\,.
\end{array}\right.
\end{aligned}\right.
\tag{\text{TVEGM}}
\end{equation*}

On the other hand, in recent years, there has been tremendous interest in developing fast iterative algorithms. Many authors have used inertial methods to devise a large number of iterative algorithms that can improve the convergence speed; see, for example,  \cite{FISTA,iFB,iDR,wtjc,TXLjnca,tanmath,liunm} and the references therein.

Motivated and inspired by the above works, in this paper, we introduce two inertial Mann-type extragradient algorithms, which are devised  with a new step size, for solving the variational inequality problem  with a monotone and Lipschitz continuous operator in real Hilbert spaces. Our algorithms  work without the knowledge of the Lipschitz constant of the involving mapping. Under some mild conditions, we prove that the iterative sequence generated by our algorithms converges   to a solution of \eqref{VIP} in norm. Some numerical experiments are provided to support the theoretical results. Our numerical results show that our new algorithms have a better convergence speed than the existing ones.

The remainder of this paper is organized as follows. In Section~\ref{sec2}, one recalls some preliminary results and lemmas for further use. Section~\ref{sec3} analyzes the convergence of the proposed algorithms. In Section~\ref{sec4}, some numerical examples are presented to illustrate the numerical behavior of the proposed algorithms and compare them with some existing ones. Finally, a brief summary is given in Section~\ref{sec5}, the last sectioin.
\section{Preliminaries}\label{sec2}
Let $C$ be a convex closed  subset of a real Hilbert space $H$. The weak convergence, which the convergence in the weak topology, and strong convergence (convergence in norm) of $\left\{x_{n}\right\}_{n=1}^{\infty}$ to $x$ are represented by $x_{n} \rightharpoonup x$ and $x_{n} \rightarrow x$, respectively. For each $x, y,z \in H$, we have the following facts:
\begin{enumerate}
\item $\|x+y\|^{2} \leq\|x\|^{2}+2\langle y, x+y\rangle$;
\item $\|\alpha x+(1-\alpha) y\|^{2}+\alpha(1-\alpha)\|x-y\|^{2}=\alpha\|x\|^{2}+(1-\alpha)\|y\|^{2}$, $\alpha \in \mathbb{R}$;
\item $\|\alpha x+\beta y+\gamma z\|^{2}= \alpha\|x\|^{2}+\beta\|y\|^{2}+\gamma\|z\|^{2}-\alpha \beta\|x-y\|^{2} -\alpha \gamma\|x-z\|^{2}-\beta \gamma\|y-z\|^{2} $, where $\alpha, \beta, \gamma \in[0, 1]$ with $\alpha+\beta+\gamma=1$.
\end{enumerate}

For every point $x \in H$, there exists a unique nearest point in $C$, denoted by $P_{C} x$ such that $P_{C}x:= \operatorname{argmin}\{\|x-y\|,\, y \in C\}$. $P_{C}$ is called the metric projection of $H$ onto $C$. It is known that $P_{C}$ is nonexpansive and $P_{C} x$ has the following basic properties:
\begin{itemize}[leftmargin=1em]
\item $ \langle x-P_{C} x, y-P_{C} x\rangle \leq 0, \, \forall y \in C $;
\item $\left\|P_{C} x-P_{C} y\right\|^{2} \leq\left\langle P_{C} x-P_{C} y, x-y\right\rangle, \,\forall y \in H$.
\end{itemize}

To prove the convergence of the proposed algorithms, we need the following lemmas.
\begin{lemma}[\cite{HSEGM}]\label{lem21}
Let $A : H \rightarrow H$ be a monotone and $ L $-Lipschitz continuous mapping on $ C $. Let $S=P_{C}(I-\mu A)$, where $\mu>0$. If $\left\{x_{n}\right\}$ is a sequence in $H$ satisfying $x_{n} \rightharpoonup q$ and $x_{n}-S x_{n} \rightarrow 0$, then $q \in \mathrm{VI}(C, A)=\operatorname{Fix}(S)$.
\end{lemma}

\begin{lemma}[\cite{Mainge2}]\label{lem22}
Assume that  $\left\{a_{n}\right\}$ is a  nonnegative real number sequence and there is a subsequence $\left\{a_{n_{j}}\right\}$ of $\left\{a_{n}\right\}$ such that $a_{n_{j}}<a_{n_{j}+1}$ for all $j \in \mathbb{N}$. Then, there exists a nondecreasing sequence $\left\{m_{k}\right\}$ of $\mathbb{N}$ such that $\lim _{k \rightarrow \infty} m_{k}=\infty$ and the following properties are satisfied by all (sufficiently large) number $k \in \mathbb{N}$ :
\[
a_{m_{k}} \leq a_{m_{k}+1} \text { and } a_{k} \leq a_{m_{k}+1}\,.
\]
In fact, $m_{k}$ is the largest number $n$ in the set $\{1,2, \ldots, k\}$ such that $a_{n}<a_{n+1}$.
\end{lemma}

\begin{lemma}[\cite{xu}]\label{lem23}
Let $\left\{a_{n}\right\}$ be  non-negative real number sequence, which satisfies
\[
a_{n+1} \leq \alpha_{n} b_{n}+\left(1-\alpha_{n}\right) a_{n}, \quad \forall n>0\,,
\]
where $\left\{\alpha_{n}\right\} \subset(0,1)$ and $\left\{b_{n}\right\}$ are a sequence such that $\sum_{n=0}^{\infty} \alpha_{n}=\infty$ and $\limsup_{n \rightarrow \infty} b_{n} \leq 0$. Then, $\lim _{n \rightarrow \infty} a_{n}=0$.
\end{lemma}

\section{Main results}\label{sec3}
In this section, we introduce two new inertial extragradient algorithms with a new step size for solving variational inequality problems and analyze their convergence. First, we assume that our proposed algorithms satisfy the following conditions.

\begin{enumerate}[label=(C\arabic*)]
\item The mapping $A: H \rightarrow H$ is monotone and $L$-Lipschitz continuous on $H$. \label{con1}
\item The solution set of the \eqref{VIP} is nonempty, that is, $\mathrm{VI}(C, A) \neq \emptyset$. \label{con2}
\item Let $ \{\epsilon_{n}\} $ be a positive sequence such that $\lim_{n \rightarrow \infty} \frac{\epsilon_{n}}{\alpha_{n}}=0$, where $ \{\alpha_{n}\}\subset (0,1) $ is with the restrictions that   $\sum_{n=1}^{\infty} \alpha_{n}=\infty$ and $\lim _{n \rightarrow \infty} \alpha_{n}=0$. Let $\left\{\beta_{n}\right\} \subset(a, b) \subset\left(0,1-\alpha_{n}\right)$ for some $a>0, b>0$. \label{con3}
\end{enumerate}

\subsection{The Mann-type inertial subgradient extragradient algorithm}
Now, we introduce a Mann-type inertial subgradient extragradient algorithm for solving variational inequality problems. The algorithm is of the form:
\begin{algorithm}[H]
\caption{The Mann-type inertial subgradient extragradient algorithm for \eqref{VIP}}
\label{alg1}
\begin{algorithmic}
	\STATE {\textbf{Initialization:} Given $ \theta>0 $, $\lambda_{1}>0$, $\mu \in(0,1)$. Let $x_{0},x_{1} \in H$ be arbitrarily fixed.}
    \STATE \textbf{Iterative Steps}: Calculate $ x_{n+1} $ as follows:
    \STATE \textbf{Step 1.} Given the iterates $x_{n-1}$ and $x_{n} (n \geq 1) $. Set
    \[w_{n}=x_{n}+\theta_{n}\left(x_{n}-x_{n-1}\right)\,,\]
    where
    \begin{equation}\label{alpha}
    \theta_{n}=\left\{\begin{array}{ll}
    \min \bigg\{\dfrac{\epsilon_{n}}{\left\|x_{n}-x_{n-1}\right\|}, \theta\bigg\}, & \text { if } x_{n} \neq x_{n-1}\,; \\
    \theta, & \text { otherwise}\,.
    \end{array}\right.
    \end{equation}
    \STATE \textbf{Step 2.} Compute
    \[y_{n}=P_{C}\left(w_{n}-\lambda_{n} A w_{n}\right)\,.\]
    If $w_{n}=y_{n}$, then stop, and $y_{n}$ is a solution of $\mathrm{VI}(C, A)$. Otherwise:
    \STATE \textbf{Step 3.} Compute
    \[z_{n}=P_{T_{n}}\left(w_{n}-\lambda_{n} A y_{n}\right)\,,\]
    where $ T_{n}:=\left\{x \in H \mid \left\langle w_{n}-\lambda_{n} A w_{n}-y_{n}, x-y_{n}\right\rangle \leq 0\right\} $.
    \STATE \textbf{Step 4.} Compute
    \[x_{n+1}=\left(1-\alpha_{n}-\beta_{n}\right) w_{n}+\beta_{n} z_{n} \,,\]
    and update
    \begin{equation}\label{lambda2}
       	\lambda_{n+1}=\left\{\begin{array}{ll}
       		\min \left\{\dfrac{\mu\left\|w_{n}-y_{n}\right\|}{\left\|A w_{n}-A y_{n}\right\|}, \lambda_{n}\right\}, & \text { if } A w_{n}-A y_{n} \neq 0\,; \\
       		\lambda_{n}, & \text { otherwise}\,.
       	\end{array}\right.
    \end{equation}
    Set $n:=n+1$ and go to \textbf{Step 1}.
\end{algorithmic}
\end{algorithm}
\begin{remark}\label{rem31}
It is easy to see from \eqref{alpha} that    
\[
\lim _{n \rightarrow \infty} \frac{\theta_{n}}{\alpha_{n}}\left\|x_{n}-x_{n-1}\right\|=0\,.
\]
Indeed, we have $\theta_{n}\left\|x_{n}-x_{n-1}\right\| \leq \epsilon_{n}$ for all $n$, which together with $\lim _{n \rightarrow \infty} \frac{\epsilon_{n}}{\alpha_{n}}=0$ implies that
\[
\lim _{n \rightarrow \infty} \frac{\theta_{n}}{\alpha_{n}}\left\|x_{n}-x_{n-1}\right\| \leq \lim _{n \rightarrow \infty} \frac{\epsilon_{n}}{\alpha_{n}}=0\,.
\]
\end{remark}
The following lemmas are quite helpful to analyze the convergence of the algorithm.
\begin{lemma}\label{lem31}
The sequence $\left\{\lambda_{n}\right\}$ generated by \eqref{lambda2} is a nonincreasing sequence and
\[
\lim _{n \rightarrow \infty} \lambda_{n}=\lambda \geq \min \Big\{\lambda_{1}, \frac{\mu}{L}\Big\}\,.
\]
\end{lemma}
\begin{proof}
It follows from \eqref{lambda2} that $\lambda_{n+1} \leq \lambda_{n}$ for all $n \in \mathbb{N} $. Hence, $\left\{\lambda_{n}\right\}$ is nonincreasing. On the other hand, we get $\left\|A w_{n}-A y_{n}\right\| \leq L\left\|w_{n}-y_{n}\right\|$ since $A$ is $L$-Lipschitz continuous. Consequently
\[
\mu \frac{\left\|w_{n}-y_{n}\right\|}{\left\|A w_{n}-A y_{n}\right\|} \geq \frac{\mu}{L} \,\,,\text {  if  }\,\, A w_{n} \neq A y_{n}\,,
\]
which together with \eqref{lambda2} implies that $ \lambda_{n} \geq \min \{\lambda_{1}, \frac{\mu}{L}\} $. Since $ \{\lambda_{n}\} $  is nonincreasing and lower bounded, we have $\lim _{n \rightarrow \infty} \lambda_{n}=\lambda \geq \min \big\{\lambda_{1}, \frac{\mu}{L}\big\}$.
\end{proof}
\begin{lemma}\label{lem32}
Assume that the Conditions \ref{con1} and \ref{con2} hold. Let $\left\{z_{n}\right\}$ be a sequence generated by Algorithm~\ref{alg1}. Then
\begin{equation}\label{q}
\left\|z_{n}-p\right\|^{2} \leq\left\|w_{n}-p\right\|^{2}-\Big(1-\mu \frac{\lambda_{n}}{\lambda_{n+1}}\Big)\left\|y_{n}-w_{n}\right\|^{2}-\Big(1-\mu \frac{\lambda_{n}}{\lambda_{n+1}}\Big)\left\|z_{n}-y_{n}\right\|^{2}
\end{equation}
for all $p \in \mathrm{V I}(C, A)$.
\end{lemma}
\begin{proof}
By the definition of  $ \{\lambda_{n}\} $, one has
\[
\left\|A w_{n}-A y_{n}\right\| \leq \frac{\mu}{\lambda_{n+1}}\left\|w_{n}-y_{n}\right\|, \quad \forall n \geq 0\,.
\]	
Using $p \in \mathrm{VI}(C, A) \subset C \subset T_{n}$, we have
\[
\begin{aligned}
2\left\|z_{n}-p\right\|^{2}=&2\left\|P_{T_{n}}\left(w_{n}-\lambda_{n} A y_{n}\right)-P_{T_{n}} p\right\|^{2} \leq2\left\langle z_{n}-p, w_{n}-\lambda_{n} A y_{n}-p\right\rangle \\
=& \left\|z_{n}-p\right\|^{2}+\left\|w_{n}-\lambda_{n} A y_{n}-p\right\|^{2}-\left\|z_{n}-w_{n}+\lambda_{n} A y_{n}\right\|^{2} \\
=& \left\|z_{n}-p\right\|^{2}+\left\|w_{n}-p\right\|^{2}+ \lambda_{n}^{2}\left\|A y_{n}\right\|^{2}-2\left\langle w_{n}-p, \lambda_{n} A y_{n}\right\rangle \\
&-\left\|z_{n}-w_{n}\right\|^{2}- \lambda_{n}^{2}\left\|A y_{n}\right\|^{2}-2\left\langle z_{n}-w_{n}, \lambda_{n} A y_{n}\right\rangle \\
=& \left\|z_{n}-p\right\|^{2}+\left\|w_{n}-p\right\|^{2}-\left\|z_{n}-w_{n}\right\|^{2}-2\left\langle z_{n}-p, \lambda_{n} A y_{n}\right\rangle\,,
\end{aligned}
\]
which implies that
\begin{equation}\label{a}
\left\|z_{n}-p\right\|^{2} \leq\left\|w_{n}-p\right\|^{2}-\left\|z_{n}-w_{n}\right\|^{2}-2\left\langle z_{n}-p, \lambda_{n} A y_{n}\right\rangle\,.
\end{equation}
We have $\left\langle A p, y_{n}-p\right\rangle \geq 0 $ since $p \in \mathrm{VI}(C, A)$. In addition, since $A$ is monotone, we have $2 \lambda_{n}\left\langle A y_{n}-A p, y_{n}-p\right\rangle \geq 0 $. Thus, adding this item to the right side of \eqref{a}, we get
\begin{equation}\label{z}
\begin{aligned}
\left\|z_{n}-p\right\|^{2} \leq &\left\|w_{n}-p\right\|^{2}-\left\|z_{n}-w_{n}\right\|^{2}-2\left\langle z_{n}-p, \lambda_{n} A y_{n}\right\rangle+2 \lambda_{n}\left\langle A y_{n}-A p, y_{n}-p\right\rangle \\
=&\left\|w_{n}-p\right\|^{2}-\left\|z_{n}-w_{n}\right\|^{2}+2\left\langle y_{n}-z_{n}, \lambda_{n} A y_{n}\right\rangle-2 \lambda_{n}\left\langle A p, y_{n}-p\right\rangle \\
\leq&\left\|w_{n}-p\right\|^{2}-\left\|z_{n}-w_{n}\right\|^{2}+2 \lambda_{n}\left\langle y_{n}-z_{n}, A y_{n}-A w_{n}\right\rangle \\
&+2 \lambda_{n}\left\langle A w_{n}, y_{n}-z_{n}\right\rangle\,.
\end{aligned}
\end{equation}
Note that 
\begin{equation}\label{w}
\begin{aligned}
2 \lambda_{n}\left\langle y_{n}-z_{n}, A y_{n}-A w_{n}\right\rangle & \leq 2 \lambda_{n}\left\|A y_{n}-A w_{n}\right\|\left\|y_{n}-z_{n}\right\| \leq 2 \mu \frac{\lambda_{n}}{\lambda_{n+1}}\left\|w_{n}-y_{n}\right\|\left\|y_{n}-z_{n}\right\| \\
& \leq \mu \frac{\lambda_{n}}{\lambda_{n+1}}\left\|w_{n}-y_{n}\right\|^{2}+\mu \frac{\lambda_{n}}{\lambda_{n+1}}\left\|y_{n}-z_{n}\right\|^{2}\,.
\end{aligned}
\end{equation}
Next, we estimate $2 \lambda_{n}\left\langle A w_{n}, y_{n}-z_{n}\right\rangle$. Since $z_{n}=P_{T_{n}}\left(w_{n}-\lambda_{n} A y_{n}\right)$ and hence $z_{n} \in T_{n}$, we have
\[
\left\langle w_{n}-\lambda_{n} A w_{n}-y_{n}, z_{n}-y_{n}\right\rangle \leq 0\,,
\]
which implies that
\begin{equation}\label{s}
\begin{aligned}
	2 \lambda_{n}\left\langle A w_{n}, y_{n}-z_{n}\right\rangle & \leq 2\left\langle y_{n}-w_{n}, z_{n}-y_{n}\right\rangle \\
	&=\left\|z_{n}-w_{n}\right\|^{2}-\left\|y_{n}-w_{n}\right\|^{2}-\left\|z_{n}-y_{n}\right\|^{2}\,.
\end{aligned}
\end{equation}
Substituting \eqref{w} and \eqref{s} into \eqref{z}, we obtain
\[
\left\|z_{n}-p\right\|^{2} \leq\left\|w_{n}-p\right\|^{2}-\Big(1-\mu \frac{\lambda_{n}}{\lambda_{n+1}}\Big)\left\|y_{n}-w_{n}\right\|^{2}-\Big(1-\mu \frac{\lambda_{n}}{\lambda_{n+1}}\Big)\left\|z_{n}-y_{n}\right\|^{2}\,.
\]
\end{proof}
\begin{theorem}\label{thm31}
Assume that Conditions \ref{con1}--\ref{con3} hold. Then the sequence $\left\{x_{n}\right\}$ generated by Algorithm~\ref{alg1} converges to $p \in  \mathrm{VI}(C, A)$  in norm, where $\|p\|=\min \{\|z\|: z \in \mathrm{VI}(C, A)\}$.
\end{theorem}
\begin{proof}
According to Lemma~\ref{lem31}, it follows that $\lim _{n \rightarrow \infty}\big(1-\mu \frac{\lambda_{n}}{\lambda_{n+1}}\big)=1-\mu>0$. Thus, there exists $n_{0} \in \mathbb{N}$  such that
\begin{equation}\label{eqq}
1-\mu \frac{\lambda_{n}}{\lambda_{n+1}}>0,\quad  \forall n \geq n_{0}\,.
\end{equation}
Combining Lemma~\ref{lem32} and \eqref{eqq}, we obtain
\begin{equation}\label{eqa}
\left\|z_{n}-p\right\| \leq\left\|w_{n}-p\right\|, \quad \forall n \geq n_{0}\,.
\end{equation}
\noindent\textbf{Claim 1.} The sequence $\left\{x_{n}\right\}$ is bounded. By the definition of $ \{x_{n+1}\} $, one has
\begin{equation}\label{eqz}
\begin{aligned}
\left\|x_{n+1}-p\right\| &=\left\|\left(1-\alpha_{n}-\beta_{n}\right) w_{n}+\beta_{n} z_{n}-p\right\| \\
&=\left\|\left(1-\alpha_{n}-\beta_{n}\right)\left(w_{n}-p\right)+\beta_{n}\left(z_{n}-p\right)-\alpha_{n} p\right\| \\
& \leq\left\|\left(1-\alpha_{n}-\beta_{n}\right)\left(w_{n}-p\right)+\beta_{n}\left(z_{n}-p\right)\right\|+\alpha_{n}\|p\|\,.
\end{aligned}
\end{equation}
On the other hand, it follows from \eqref{eqa} that
\[
\begin{aligned}
&\quad\left\|\left(1-\alpha_{n}-\beta_{n}\right)\left(w_{n}-p\right)+\beta_{n}\left(z_{n}-p\right)\right\|^{2} \\
&=\left(1-\alpha_{n}-\beta_{n}\right)^{2}\left\|w_{n}-p\right\|^{2}+2\left(1-\alpha_{n}-\beta_{n}\right) \beta_{n}\left\langle w_{n}-p, z_{n}-p\right\rangle+\beta_{n}^{2}\left\|z_{n}-p\right\|^{2} \\
& \leq\left(1-\alpha_{n}-\beta_{n}\right)^{2}\left\|w_{n}-p\right\|^{2}+2\left(1-\alpha_{n}-\beta_{n}\right) \beta_{n}\left\|z_{n}-p\right\|\left\|w_{n}-p\right\|+\beta_{n}^{2}\left\|z_{n}-p\right\|^{2} \\
& \leq\left(1-\alpha_{n}-\beta_{n}\right)^{2}\left\|w_{n}-p\right\|^{2}+2\left(1-\alpha_{n}-\beta_{n}\right) \beta_{n}\left\|w_{n}-p\right\|^{2}+\beta_{n}^{2}\left\|w_{n}-p\right\|^{2}  \\
&=\left(1-\alpha_{n}\right)^{2}\left\|w_{n}-p\right\|^{2},\quad \forall n \geq n_{0}\,,
\end{aligned}
\]
which yields
\begin{equation}\label{eqw}
\left\|\left(1-\alpha_{n}-\beta_{n}\right)\left(w_{n}-p\right)+\beta_{n}\left(z_{n}-p\right)\right\| \leq\left(1-\alpha_{n}\right)\left\|w_{n}-p\right\|,\,\, \forall n \geq n_{0}\,.
\end{equation}
Using the definition of $ w_{n} $, we can write
\begin{equation}\label{c}
\begin{aligned}
\left\|w_{n}-p\right\| 
&\leq\left\|x_{n}-p\right\|+\alpha_{n} \cdot \frac{\theta_{n}}{\alpha_{n}}\left\|x_{n}-x_{n-1}\right\|\,.
\end{aligned}
\end{equation}
By Remark~\ref{rem31}, we have $\frac{\theta_{n}}{\alpha_{n}}\left\|x_{n}-x_{n-1}\right\| \rightarrow 0$. Thus, there exists a constant $M_{1}>0$ such that
\begin{equation}\label{r}
\frac{\theta_{n}}{\alpha_{n}}\left\|x_{n}-x_{n-1}\right\| \leq M_{1}, \quad \forall n \geq 1\,.
\end{equation}
From \eqref{eqa}, \eqref{c} and \eqref{r}, we find that
\begin{equation}\label{f}
\left\|z_{n}-p\right\| \leq\left\|w_{n}-p\right\| \leq\left\|x_{n}-p\right\|+\alpha_{n} M_{1},\quad \forall n \geq n_{0}\,.
\end{equation}
Combining \eqref{eqz}, \eqref{eqw} and \eqref{f}, we deduce that
\begin{equation*}
\begin{aligned}
\left\|x_{n+1}-p\right\| & \leq\left(1-\alpha_{n}\right)\left\|w_{n}-p\right\|+\alpha_{n}\|p\| \\
& \leq \left(1-\alpha_{n}\right)\left\|x_{n}-p\right\|+\alpha_{n}(\|p\|+M_{1})\\
& \leq \max \left\{\left\|x_{n}-p\right\|,\|p\|+M_{1}\right\} \\
& \leq \cdots \leq \max \left\{\left\|x_{n_{0}}-p\right\|,\|p\|+M_{1}\right\}\,.
\end{aligned}
\end{equation*}	
That is, the sequence $\left\{x_{n}\right\}$ is bounded. So the sequences $\left\{w_{n}\right\}$ and $ \{z_{n}\} $ are also bounded.

\noindent\textbf{Claim 2.}
\[
\begin{aligned}
&\quad\beta_{n} \Big(1-\mu \frac{\lambda_{n}}{\lambda_{n+1}}\Big)\left\|w_{n}-y_{n}\right\|^{2}+\beta_{n}\Big(1-\mu \frac{\lambda_{n}}{\lambda_{n+1}}\Big)\left\|y_{n}-z_{n}\right\|^{2} \\
&\leq\left\|x_{n}-p\right\|^{2}-\left\|x_{n+1}-p\right\|^{2}+\alpha_{n}(\|p\|^{2}+M_{2})\,.
\end{aligned}
\]
Indeed, by the definition of $ x_{n+1} $, one obtains
\begin{equation}\label{eqs}
\begin{aligned}
\left\|x_{n+1}-p\right\|^{2}=&\left\|\left(1-\alpha_{n}-\beta_{n}\right) w_{n}+\beta_{n} z_{n}-p\right\|^{2} \\
=&\left\|\left(1-\alpha_{n}-\beta_{n}\right)\left(w_{n}-p\right)+\beta_{n}\left(z_{n}-p\right)+\alpha_{n}(-p)\right\|^{2} \\
=&\left(1-\alpha_{n}-\beta_{n}\right)\left\|w_{n}-p\right\|^{2}+\beta_{n}\left\|z_{n}-p\right\|^{2}+\alpha_{n}\|p\|^{2}\\
&-\beta_{n}\left(1-\alpha_{n}-\beta_{n}\right)\left\|w_{n}-z_{n}\right\|^{2}-\alpha_{n}\left(1-\alpha_{n}-\beta_{n}\right)\left\|w_{n}\right\|^{2}-\alpha_{n} \beta_{n}\left\|z_{n}\right\|^{2} \\
\leq&\left(1-\alpha_{n}-\beta_{n}\right)\left\|w_{n}-p\right\|^{2}+\beta_{n}\left\|z_{n}-p\right\|^{2}+\alpha_{n}\|p\|^{2}\,.
\end{aligned}
\end{equation}
In view of \eqref{f}, one sees that
\begin{equation}\label{eqo}
\begin{aligned}
\left\|w_{n}-p\right\|^{2} & \leq\left(\left\|x_{n}-p\right\|+\alpha_{n} M_{1}\right)^{2} \\
&=\left\|x_{n}-p\right\|^{2}+\alpha_{n}\left(2 M_{1}\left\|x_{n}-p\right\|+\alpha_{n} M_{1}^{2}\right) \\
& \leq\left\|x_{n}-p\right\|^{2}+\alpha_{n} M_{2}
\end{aligned}
\end{equation}
for some $M_{2}>0 $. Thus, using Lemma~\ref{lem32}, \eqref{eqs} and \eqref{eqo}, we obtain
\[
\begin{aligned}
\left\|x_{n+1}-p\right\|^{2} \leq &\left(1-\alpha_{n}-\beta_{n}\right)\left\|w_{n}-p\right\|^{2}+\beta_{n}\left\|w_{n}-p\right\|^{2}-\beta_{n}\Big(1-\mu \frac{\lambda_{n}}{\lambda_{n+1}}\Big)\left\|w_{n}-y_{n}\right\|^{2} \\
&-\beta_{n}\Big(1-\mu \frac{\lambda_{n}}{\lambda_{n+1}}\Big)\left\|y_{n}-z_{n}\right\|^{2}+\alpha_{n}\|p\|^{2} \\
\leq&\left\|x_{n}-p\right\|^{2}-\beta_{n}\Big(1-\mu\frac{\lambda_{n}}{\lambda_{n+1}}\Big)\left\|w_{n}-y_{n}\right\|^{2} \\
&-\beta_{n}\Big(1-\mu \frac{\lambda_{n}}{\lambda_{n+1}}\Big)\left\|y_{n}-z_{n}\right\|^{2}+\alpha_{n}(\|p\|^{2}+M_{2})\,.
\end{aligned}
\]
\noindent\textbf{Claim 3.}
\[
\begin{aligned}
\left\|x_{n+1}-p\right\|^{2} \leq&\left(1-\alpha_{n}\right)\left\|x_{n}-p\right\|^{2}+\alpha_{n}\Big[2 \beta_{n}\left\|w_{n}-z_{n}\right\|\left\|x_{n+1}-p\right\| \Big.\\
&\Big.+2\left\langle p, p-x_{n+1}\right\rangle + \frac{3 M \theta_{n}}{ \alpha_{n}}\left\|x_{n}-x_{n-1}\right\|\Big],\,\, \forall n \geq n_{0}\,.
\end{aligned}
\]
Indeed, by the definition of $ w_{n} $, one obtains
\begin{equation}\label{eqk}
\begin{aligned}
\left\|w_{n}-p\right\|^{2} &=\left\|x_{n}+\theta_{n}\left(x_{n}-x_{n-1}\right)-p\right\|^{2} \\
&=\left\|x_{n}-p\right\|^{2}+2 \theta_{n}\left\langle x_{n}-p, x_{n}-x_{n-1}\right\rangle+\theta_{n}^{2}\left\|x_{n}-x_{n-1}\right\|^{2} \\
&\leq\left\|x_{n}-p\right\|^{2} + 3M\theta_{n}\left\|x_{n}-x_{n-1}\right\|\,,
\end{aligned}
\end{equation}
where $M:=\sup _{n \in \mathbb{N}}\left\{\left\|x_{n}-p\right\|, \theta\left\|x_{n}-x_{n-1}\right\|\right\}>0$. Setting $t_{n}=\left(1-\beta_{n}\right) w_{n}+\beta_{n} z_{n}$, one has
\begin{equation}\label{eqx}
\left\|t_{n}-w_{n}\right\|=\beta_{n}\left\|w_{n}-z_{n}\right\|\,.
\end{equation}
It follows from \eqref{f} that
\begin{equation}\label{eqe}
\begin{aligned}
\left\|t_{n}-p\right\| &=\left\|\left(1-\beta_{n}\right)\left(w_{n}-p\right)+\beta_{n}\left(z_{n}-p\right)\right\| \\
&\leq\left(1-\beta_{n}\right)\left\|w_{n}-p\right\|+\beta_{n}\left\|w_{n}-p\right\|  \\
&=\left\|w_{n}-p\right\|, \quad \forall n \geq n_{0}\,.
\end{aligned}
\end{equation}
From \eqref{eqk}, \eqref{eqx} and \eqref{eqe}, for all $ n \geq n_{0} $, we get
\[
\begin{aligned}
\left\|x_{n+1}-p\right\|^{2}
=&\left\|\left(1-\beta_{n}\right) w_{n}+\beta_{n} z_{n}-\alpha_{n} w_{n}-p\right\|^{2} \\
=&\left\|\left(1-\alpha_{n}\right)\left(t_{n}-p\right)-\alpha_{n}\left(w_{n}-t_{n}\right)-\alpha_{n} p\right\|^{2} \\
\leq&\left(1-\alpha_{n}\right)^{2}\left\|t_{n}-p\right\|^{2}-2\alpha_{n}\left\langle w_{n}-t_{n}+ p, x_{n+1}-p\right\rangle \\
=&\left(1-\alpha_{n}\right)^{2}\left\|t_{n}-p\right\|^{2}+2 \alpha_{n}\left\langle w_{n}-t_{n}, p-x_{n+1}\right\rangle+2 \alpha_{n}\left\langle p, p-x_{n+1}\right\rangle \\
\leq&\left(1-\alpha_{n}\right)\left\|t_{n}-p\right\|^{2}+2 \alpha_{n}\left\|w_{n}-t_{n}\right\|\left\|x_{n+1}-p\right\|+2 \alpha_{n}\left\langle p, p-x_{n+1}\right\rangle \\
\leq&\left(1-\alpha_{n}\right)\left\|x_{n}-p\right\|^{2}+\alpha_{n}\Big[2 \beta_{n}\left\|w_{n}-z_{n}\right\|\left\|x_{n+1}-p\right\| \Big.\\
&\Big.+2\left\langle p, p-x_{n+1}\right\rangle + \frac{3 M \theta_{n}}{ \alpha_{n}}\left\|x_{n}-x_{n-1}\right\|\Big]\,.
\end{aligned}
\]
\noindent\textbf{Claim 4.} The sequence $\{\left\|x_{n}-p\right\|^{2}\}$ converges to zero by considering two possible cases on the sequence $\{\left\|x_{n}-p\right\|^{2}\}$.

\textbf{Case 1.} There exists an $N \in \mathbb{N}$, such that $\|x_{n+1}-p\|^{2} \leq\|x_{n}-p\|^{2}$ for all $n \geq N$. This implies that $\lim _{n \rightarrow \infty}\left\|x_{n}-p\right\|^{2}$ exists. In view of  $\lim _{n \rightarrow \infty}\big(1-\mu \frac{\lambda_{n}}{\lambda_{n}+1}\big)=1-\mu>0$ and Condition \ref{con3}. It implies from Claim 2 that
\[
\lim _{n \rightarrow \infty}\left\|w_{n}-y_{n}\right\|=0, \text{  and  } \lim _{n \rightarrow \infty}\left\|y_{n}-z_{n}\right\|=0\,.
\]
This implies that $ \lim _{n \rightarrow \infty}\left\|z_{n}-w_{n}\right\|=0 $, which, together with the boundedness of $ \{x_{n}\} $, it is further concluded that
\[
\lim _{n \rightarrow \infty} \beta_{n}\left\|w_{n}-z_{n}\right\|\left\|x_{n+1}-p\right\|=0\,.
\]
According to the definition of $ w_{n} $, one has
\begin{equation*}\label{eqr}
\left\|x_{n}-w_{n}\right\|=\theta_{n}\left\|x_{n}-x_{n-1}\right\|=\alpha_{n} \cdot \frac{\theta_{n}}{\alpha_{n}}\left\|x_{n}-x_{n-1}\right\| \rightarrow 0 \text { as } n \rightarrow \infty\,.
\end{equation*}
On the other hand, one sees that
\[
\left\|x_{n+1}-w_{n}\right\| \leq \alpha_{n}\left\|w_{n}\right\|+\beta_{n}\left\|z_{n}-w_{n}\right\| \rightarrow 0 \text { as } n \rightarrow \infty\,.
\]
This together with $ \lim _{n \rightarrow \infty}\left\|x_{n}-w_{n}\right\|=0 $ implies that
\begin{equation*}\label{y}
\lim _{n \rightarrow \infty}\left\|x_{n+1}-x_{n}\right\|=0\,.
\end{equation*}
Since $\left\{x_{n}\right\}$ is bounded, there exists a subsequence $\{x_{n_{j}}\}$ of $\left\{x_{n}\right\}$, such that $x_{n_{j}} \rightharpoonup q$ and
\[
\limsup _{n \rightarrow \infty}\left\langle p, p-x_{n}\right\rangle=\lim _{j \rightarrow \infty}\left\langle p, p-x_{n_{j}}\right\rangle=\langle p, p-q\rangle\,.
\]
We get $w_{n_{j}} \rightharpoonup q$ since $ \left\|x_{n}-w_{n}\right\|\rightarrow 0 $, this together with $\lim _{n \rightarrow \infty} \lambda_{n}=\lambda>0 $ and $\left\|w_{n}-y_{n}\right\|\rightarrow 0$, in the light of Lemma~\ref{lem21}, we obtain $q \in \mathrm{VI}(C, A)$. Since $q \in \operatorname{VI}(C, A)$ and $\|p\|=\min \{\|z\|: z \in \operatorname{VI}(C, A)\}$, that is $p=P_{\mathrm{VI}(C, A)} 0$, we
deduce that
\[
\limsup _{n \rightarrow \infty}\left\langle p, p-x_{n}\right\rangle=\langle p, p-q\rangle \leq 0\,.
\]
From $\left\|x_{n+1}-x_{n}\right\| \rightarrow 0$, we get
\[
\limsup _{n \rightarrow \infty}\left\langle p, p-x_{n+1}\right\rangle \leq 0\,.
\]
Therefore, using Claim 3 and Remark~\ref{rem31} in Lemma~\ref{lem23}, we conclude that $x_{n} \rightarrow p$.

\textbf{Case 2.} There exists a subsequence $\{\|x_{n_{j}}-p\|^{2}\}$ of $\{\|x_{n}-p\|^{2}\}$ such that $\|x_{n_{j}}-p\|^{2}<$ $\|x_{n_{j}+1}-p\|^{2}$ for all $j \in \mathbb{N}$. In this case, it follows from Lemma~\ref{lem22} that there exists a nondecreasing sequence $\left\{m_{k}\right\}$ of $\mathbb{N}$ such that $\lim _{k \rightarrow \infty} m_{k}=\infty$ and the following inequalities hold for all $k \in \mathbb{N}$ :
\[
\left\|x_{m_{k}}-p\right\|^{2} \leq\|x_{m_{k}+1}-p\|^{2}, \text { and }\left\|x_{k}-p\right\|^{2} \leq\|x_{m_{k}+1}-p\|^{2}\,.
\]
By Claim 2, we have
\[
\begin{aligned}
&\quad\beta_{m_{k}}\Big(1-\mu \frac{\lambda_{m_{k}}}{\lambda_{m_{k}+1}}\Big)\left\|w_{m_{k}}-y_{m_{k}}\right\|^{2}+\beta_{m_{k}}\Big(1-\mu \frac{\lambda_{m_{k}}}{\lambda_{m_{k}+1}}\Big)\left\|y_{m_{k}}-z_{m_{k}}\right\|^{2} \\
&\leq\|x_{m_{k}}-p\|^{2}-\|x_{m_{k}+1}-p\|^{2}+\alpha_{m_{k}}(\|p\|^{2}+M_{2}) \\
&\leq \alpha_{m_{k}}(\|p\|^{2}+M_{2})\,.
\end{aligned}
\]
Therefore, from condition \ref{con3}, we get
\[
\lim _{k \rightarrow \infty}\left\|w_{m_{k}}-y_{m_{k}}\right\|=0, \text { and } \lim _{k \rightarrow \infty}\left\|y_{m_{k}}-z_{m_{k}}\right\|=0\,.
\]
As proved in the first case, we get $ \|x_{m_{k}+1}-x_{m_{k}}\| \rightarrow 0 $  and $ \limsup_{k \rightarrow \infty} \langle p, p-x_{m_{k}+1}\rangle \leq 0 $. Since Claim 3 and $ \|x_{m_{k}}-p\|^{2} \leq\|x_{m_{k}+1}-p\|^{2} $, we have
\[
\begin{aligned}
\|x_{m_{k}+1}-p\|^{2} \leq &\left(1-\alpha_{m_{k}}\right)\|x_{m_{k}+1}-p\|^{2} +\alpha_{m_{k}}\Big[2 \beta_{m_{k}}\left\|w_{m_{k}}-z_{m_{k}}\right\|\|x_{m_{k}+1}-p\|\Big.\\
&\Big.+2\langle p, p-x_{m_{k}+1}\rangle +\frac{3 M \theta_{m_{k}}}{ \alpha_{m_{k}}}\|x_{m_{k}}-x_{m_{k}-1}\|\Big] \,.
\end{aligned}
\]
This implies that
\[
\left\|x_{k}-p\right\|^{2} \leq 2 \beta_{m_{k}}\left\|w_{m_{k}}-z_{m_{k}}\right\|\|x_{m_{k}+1}-p\|+2\left\langle p, p-x_{m_{k}+1}\right\rangle+\frac{3 M \theta_{m_{k}}}{ \alpha_{m_{k}}}\|x_{m_{k}}-x_{m_{k}-1}\|\,.
\]
Therefore, we obtain $\lim \sup _{k \rightarrow \infty}\left\|x_{k}-p\right\| \leq 0 $, that is, $x_{k} \rightarrow p$. The proof is completed.
\end{proof}
\subsection{The Mann-type inertial Tseng's extragradient algorithm}
In this section, we introduce a Mann-type inertial Tseng's extragradient algorithm for solving variational inequality problems. Our algorithm is as follows:
\begin{algorithm}[H]
\caption{The Mann-type inertial Tseng's extragradient algorithm for \eqref{VIP}}
\label{alg2}
\begin{algorithmic}
	\STATE {\textbf{Initialization:} Given $ \theta>0 $, $\lambda_{1}>0$, $\mu \in(0,1)$. Let $x_{0},x_{1} \in H$ be arbitrary.}
	\STATE \textbf{Iterative Steps}: Calculate $ x_{n+1} $ as follows:
	\STATE \textbf{Step 1.} Given the iterates $x_{n-1}$ and $x_{n}(n \geq 1) $. Set
	\[w_{n}=x_{n}+\theta_{n}\left(x_{n}-x_{n-1}\right)\,,\]
	where
	\begin{equation*}\label{alpha1}
	\theta_{n}=\left\{\begin{array}{ll}
	\min \bigg\{\dfrac{\epsilon_{n}}{\left\|x_{n}-x_{n-1}\right\|}, \theta\bigg\}, & \text { if } x_{n} \neq x_{n-1}\,; \\
	\theta, & \text { otherwise}\,.
	\end{array}\right.
	\end{equation*}
	\STATE \textbf{Step 2.} Compute
	\[y_{n}=P_{C}\left(w_{n}-\lambda_{n} A w_{n}\right)\,.\]
	If $w_{n}=y_{n}$, then stop, and $y_{n}$ is a solution of $\mathrm{VI}(C, A)$. Otherwise:
	\STATE \textbf{Step 3.} Compute
	\[z_{n}=y_{n}-\lambda_{n}\left(A y_{n}-A w_{n}\right)\,,\]
	\STATE \textbf{Step 4.} Compute
	\[x_{n+1}=\left(1-\alpha_{n}-\beta_{n}\right) w_{n}+\beta_{n} z_{n} \,,\]
	and update
	\begin{equation*}
	\lambda_{n+1}=\left\{\begin{array}{ll}
	\min \left\{\dfrac{\mu\left\|w_{n}-y_{n}\right\|}{\left\|A w_{n}-A y_{n}\right\|}, \lambda_{n}\right\}, & \text { if } A w_{n}-A y_{n} \neq 0\,; \\
	\lambda_{n}, & \text { otherwise}\,.
	\end{array}\right.
	\end{equation*}
	Set $n:=n+1$ and go to \textbf{Step 1}.
\end{algorithmic}
\end{algorithm}

The following lemma is very helpful for analyzing the convergence of the Algorithm~\ref{alg2}.
\begin{lemma}\label{lem41}
Assume that Conditions \ref{con1} and \ref{con2} hold. Let $\left\{z_{n}\right\}$ be a sequence generated by Algorithm~\ref{alg2}. Then
\begin{equation*}
\left\|z_{n}-p\right\|^{2} \leq\left\|w_{n}-p\right\|^{2}-\Big(1-\mu^{2} \frac{\lambda_{n}^{2}}{\lambda_{n+1}^{2}}\Big)\left\|w_{n}-y_{n}\right\|^{2},\quad \forall p \in \mathrm{VI}(C, A)\,,
\end{equation*}
and
\begin{equation*}
\left\|z_{n}-y_{n}\right\| \leq \mu \frac{\lambda_{n}}{\lambda_{n+1}}\left\|w_{n}-y_{n}\right\| \,.
\end{equation*}
\end{lemma}
\begin{proof}
First, using the definition of $\left\{\lambda_{n}\right\}$, it is easy to see that
\begin{equation}\label{qq}
\left\|A w_{n}-A y_{n}\right\| \leq \frac{\mu}{\lambda_{n+1}}\left\|w_{n}-y_{n}\right\|, \quad \forall n\geq 0\,.
\end{equation}
By the definition of $ z_{n} $, one sees that
\begin{equation}\label{aa}
\begin{aligned}
\left\|z_{n}-p\right\|^{2}=&\left\|y_{n}-\lambda_{n}\left(A y_{n}-A w_{n}\right)-p\right\|^{2} \\
=&\left\|w_{n}-p\right\|^{2}+\left\|y_{n}-w_{n}\right\|^{2}+2\left\langle y_{n}-w_{n}, w_{n}-p\right\rangle \\
&+\lambda_{n}^{2}\left\|A y_{n}-A w_{n}\right\|^{2}-2 \lambda_{n}\left\langle y_{n}-p, A y_{n}-A w_{n}\right\rangle \\
=&\left\|w_{n}-p\right\|^{2}+\left\|y_{n}-w_{n}\right\|^{2}-2\left\langle y_{n}-w_{n}, y_{n}-w_{n}\right\rangle+2\left\langle y_{n}-w_{n}, y_{n}-p\right\rangle \\
&+\lambda_{n}^{2}\left\|A y_{n}-A w_{n}\right\|^{2}-2 \lambda_{n}\left\langle y_{n}-p, A y_{n}-A w_{n}\right\rangle \\
=&\left\|w_{n}-p\right\|^{2}-\left\|y_{n}-w_{n}\right\|^{2}+2\left\langle y_{n}-w_{n}, y_{n}-p\right\rangle \\
&+\lambda_{n}^{2}\left\|A y_{n}-A w_{n}\right\|^{2}-2 \lambda_{n}\left\langle y_{n}-p, A y_{n}-A w_{n}\right\rangle\,.
\end{aligned}
\end{equation}
Since $y_{n}=P_{C}\left(w_{n}-\lambda_{n} A w_{n}\right)$, using the property of projection, we obtain
\[
\left\langle y_{n}-w_{n}+\lambda_{n} A w_{n}, y_{n}-p\right\rangle \leq 0\,,
\]
or equivalently
\begin{equation}\label{zz}
\left\langle y_{n}-w_{n}, y_{n}-p\right\rangle \leq-\lambda_{n}\left\langle A w_{n}, y_{n}-p\right\rangle\,.
\end{equation}
From \eqref{qq}, \eqref{aa} and \eqref{zz}, we have
\begin{equation}\label{ww}
\begin{aligned}
\left\|z_{n}-p\right\|^{2}  \leq&\left\|w_{n}-p\right\|^{2}-\left\|y_{n}-w_{n}\right\|^{2}-2 \lambda_{n}\left\langle A w_{n}, y_{n}-p\right\rangle+\mu^{2} \frac{\lambda_{n}^{2}}{\lambda_{n+1}^{2}}\left\|w_{n}-y_{n}\right\|^{2} \\
&-2 \lambda_{n}\left\langle y_{n}-p, A y_{n}-A w_{n}\right\rangle \\
=&\left\|w_{n}-p\right\|^{2}-\Big(1-\mu^{2} \frac{\lambda_{n}^{2}}{\lambda_{n+1}^{2}}\Big)\left\|w_{n}-y_{n}\right\|^{2}-2 \lambda_{n}\left\langle y_{n}-p, A y_{n}-Ap\right\rangle\\
&-2 \lambda_{n}\left\langle y_{n}-p, Ap\right\rangle\,.
\end{aligned}
\end{equation}
Since $p \in \mathrm{VI}(C, A)$ and the monotonicity of $A$, we get
\begin{equation}\label{ss}
\left\langle A p, y_{n}-p\right\rangle \geq 0 \text{  and  } \left\langle A y_{n}-Ap, y_{n}-p\right\rangle \geq 0\,.
\end{equation}
Combining \eqref{ww} and \eqref{ss}, we deduce that
\[
\left\|z_{n}-p\right\|^{2} \leq\left\|w_{n}-p\right\|^{2}-\Big(1-\mu^{2} \frac{\lambda_{n}^{2}}{\lambda_{n+1}^{2}}\Big)\left\|w_{n}-y_{n}\right\|^{2}\,.
\]
From the definition of $ z_{n} $ and \eqref{qq}, we obtain
\begin{equation*}
\left\|z_{n}-y_{n}\right\| \leq \mu \frac{\lambda_{n}}{\lambda_{n+1}}\left\|w_{n}-y_{n}\right\| \,.
\end{equation*}
\end{proof}
\begin{theorem}\label{thm41}
Assume that Conditions \ref{con1}--\ref{con3} hold. Then the sequence $\left\{x_{n}\right\}$ generated by Algorithm~\ref{alg2} converges to  $p \in  \mathrm{VI}(C, A)$  in norm, where $\|p\|=\min \{\|z\|: z \in \mathrm{VI}(C, A)\}$.
\end{theorem}
\begin{proof}
Since $\lim _{n \rightarrow \infty}\big(1-\mu^{2} \frac{\lambda_{n}^{2}}{\lambda_{n+1}^{2}}\big)=1-\mu^{2}>0$, there exists $n_{0} \in \mathbb{N}$ such that
\begin{equation}\label{eqd}
1-\mu^{2} \frac{\lambda_{n}^{2}}{\lambda_{n+1}^{2}}>0,\quad \forall n \geq n_{0}\,.
\end{equation}
Combining Lemma \ref{lem41} and \eqref{eqd}, we get
\begin{equation}
\left\|z_{n}-p\right\| \leq\left\|w_{n}-p\right\|, \quad \forall n \geq n_{0}\,.
\end{equation}
\noindent\textbf{Claim 1.} The sequence $\left\{x_{n}\right\}$ is bounded. Using the same arguments with the Claim 1 in the Theorem~\ref{thm31}, we get that $\left\{x_{n}\right\}$ is bounded. Consequently, $ \{w_{n}\} $ and $\left\{z_{n}\right\}$ are also bounded.

\noindent\textbf{Claim 2.}
\[
\begin{aligned}
&\quad\beta_{n}\Big(1-\mu^{2} \frac{\lambda_{n}^{2}}{\lambda_{n+1}^{2}}\Big)\left\|w_{n}-y_{n}\right\|^{2}+\beta_{n}\left(1-\alpha_{n}-\beta_{n}\right)\left\|w_{n}-z_{n}\right\|^{2} \\
&\leq\left\|x_{n}-p\right\|^{2}-\left\|x_{n+1}-p\right\|^{2}+\alpha_{n}(\|p\|^{2}+M_{2})\,.
\end{aligned}
\]
Indeed, by the definition of $ x_{n+1} $, we have
\begin{equation}\label{eqs1}
\begin{aligned}
\left\|x_{n+1}-p\right\|^{2}=&\left\|\left(1-\alpha_{n}-\beta_{n}\right) w_{n}+\beta_{n} z_{n}-p\right\|^{2} \\
\leq&\left(1-\alpha_{n}-\beta_{n}\right)\left\|w_{n}-p\right\|^{2}+\beta_{n}\left\|z_{n}-p\right\|^{2}+\alpha_{n}\|p\|^{2}\\
&-\beta_{n}\left(1-\alpha_{n}-\beta_{n}\right)\left\|w_{n}-z_{n}\right\|^{2}\,.
\end{aligned}
\end{equation}
Combining \eqref{eqo}, Lemma \ref{lem41} and \eqref{eqs1}, we obtain
\[
\begin{aligned}
\left\|x_{n+1}-p\right\|^{2} \leq &\left(1-\alpha_{n}-\beta_{n}\right)\left\|w_{n}-p\right\|^{2}+\beta_{n}\left\|w_{n}-p\right\|^{2}-\beta_{n}\Big(1-\mu^{2} \frac{\lambda_{n}^{2}}{\lambda_{n+1}^{2}}\Big)\left\|w_{n}-y_{n}\right\|^{2} \\
&+\alpha_{n}\|p\|^{2}-\beta_{n}\left(1-\alpha_{n}-\beta_{n}\right)\left\|w_{n}-z_{n}\right\|^{2}\\
\leq&\left\|x_{n}-p\right\|^{2}-\beta_{n}\Big(1-\mu^{2}\frac{\lambda_{n}^{2}}{\lambda_{n+1}^{2}}\Big)\left\|w_{n}-y_{n}\right\|^{2}+\alpha_{n}(\|p\|^{2}+M_{2}) \\
&-\beta_{n}\left(1-\alpha_{n}-\beta_{n}\right)\left\|w_{n}-z_{n}\right\|^{2}\,.
\end{aligned}
\]
The desired result can be obtained by a simple deformation.

\noindent\textbf{Claim 3.}
\[
\begin{aligned}
\left\|x_{n+1}-p\right\|^{2} \leq&\left(1-\alpha_{n}\right)\left\|x_{n}-p\right\|^{2}+\alpha_{n}\Big[2 \beta_{n}\left\|w_{n}-z_{n}\right\|\left\|x_{n+1}-p\right\| \Big.\\
&\Big.+2\left\langle p, p-x_{n+1}\right\rangle + \frac{3 M \theta_{n}}{ \alpha_{n}}\left\|x_{n}-x_{n-1}\right\|\Big],\,\, \forall n \geq n_{0}\,.
\end{aligned}
\]
The desired result can be obtained by using the same arguments as in the Theorem~\ref{thm31} of Claim~3.

\noindent\textbf{Claim 4.} The sequence $\{\|x_{n}-q\|^{2}\}$ converges to zero. The proof is similar to the Claim 4 in Theorem~\ref{thm31}, we leave it for the reader to verify.
\end{proof}
\section{Numerical examples}\label{sec4}
In this section, we provide some numerical examples to show the numerical behavior of our proposed algorithms, namely Algorithm~\ref{alg1} (shortly, MiSEGM) and  Algorithm~\ref{alg2} (MiTEGM), and also to compare them with some existing ones  including the Halpern subgradient extragradient algorithm \eqref{HSEGM} \cite{HSEGM}, the viscosity subgradient extragradient algorithm \eqref{VSEGM} \cite{shefu}, the Tseng's viscosity extragradient algorithm \eqref{TVEGM} \cite{TVEGM}, the Mann-type subgradient extragradient algorithm (MaSEGM) \cite{TVCAM} and the Mann-type Tseng's extragradient algorithm (MaTEGM) \cite{TVCAM}. We use the FOM Solver~\cite{FOM} to effectively calculate the projections onto $ C $ and $ T_{n} $. All the programs were implemented in MATLAB 2018a on a Intel(R) Core(TM) i5-8250U CPU @ 1.60GHz computer with RAM 8.00 GB.

Our parameters are set as follows. In all algorithms, set $ \alpha_{n}=1/(n+1) $ and $ \beta_{n}=0.5(1-\alpha_{n}) $. For the proposed algorithms and the algorithms (MaSEGM) and (MaTEGM), we choose $ \lambda_{1}=1 $, $ \mu=0.5 $. Take $ \theta=0.4 $, $ \epsilon_{n}=100/(n+1)^2 $ in our proposed algorithms. For the algorithm \eqref{VSEGM}, we choose $ \ell=0.5 $, $ \mu=0.4 $ and $ f(x)=0.9x $. Setting $ \lambda_{0}=1 $, $  \mu=0.5 $ and $ f(x)=0.9x $ in the algorithm \eqref{TVEGM}. For the algorithm \eqref{HSEGM}, we choose the step size as $ \lambda_{n}= 0.99/L$. Maximum iteration $ 200 $ as a common stopping criterion. In our numerical examples, the solution $ x^{*} $ of the problems are known, so we use $ D_{n}=\|x_{n}-x^{*}\| $ to measure the $ n $-th iteration error.
\begin{example}\label{ex1}
Let us consider the following nonlinear optimization problem via
\begin{equation}\label{opt}
\begin{array}{l}{\min \;\;F(x)=1+x_{1}^{2}-e^ {-x_{2}^{2}}} \\ {\text { s.t. }-5e \leq x \leq 5e},\end{array}
\end{equation}
where $x=\left(x_{1}, x_{2}\right)^{\mathsf{T}} \in {R}^{2}$, $ e=(1,1)^{\mathsf{T}} $.  Observe that  $\nabla F(x)=(2 x_{1}, 2x_{2}e^ {-x_{2}^{2}})^{\mathsf{T}}$ and the optimal solution for $ F(x) $ is $ x^*=(0,0)^{\mathsf{T}} $. Taking $A(x)=\nabla F(x)$, it is easy to check that  $ A(x) $ is monotone and Lipschizt continuous with  constant $L=2$  on the closed and convex subset $C=\left\{x \in {R}^{2}:-5e \leq x \leq 5e\right\}$. The initial values $ x_{0} = x_{1} $ are randomly generated by \emph{rand(2,1)} in MATLAB. The numerical results are reported in Figs.~\ref{ex1_fig1} and \ref{ex1_fig2}.
\begin{figure}[htbp]
\centering
\includegraphics[scale=0.66]{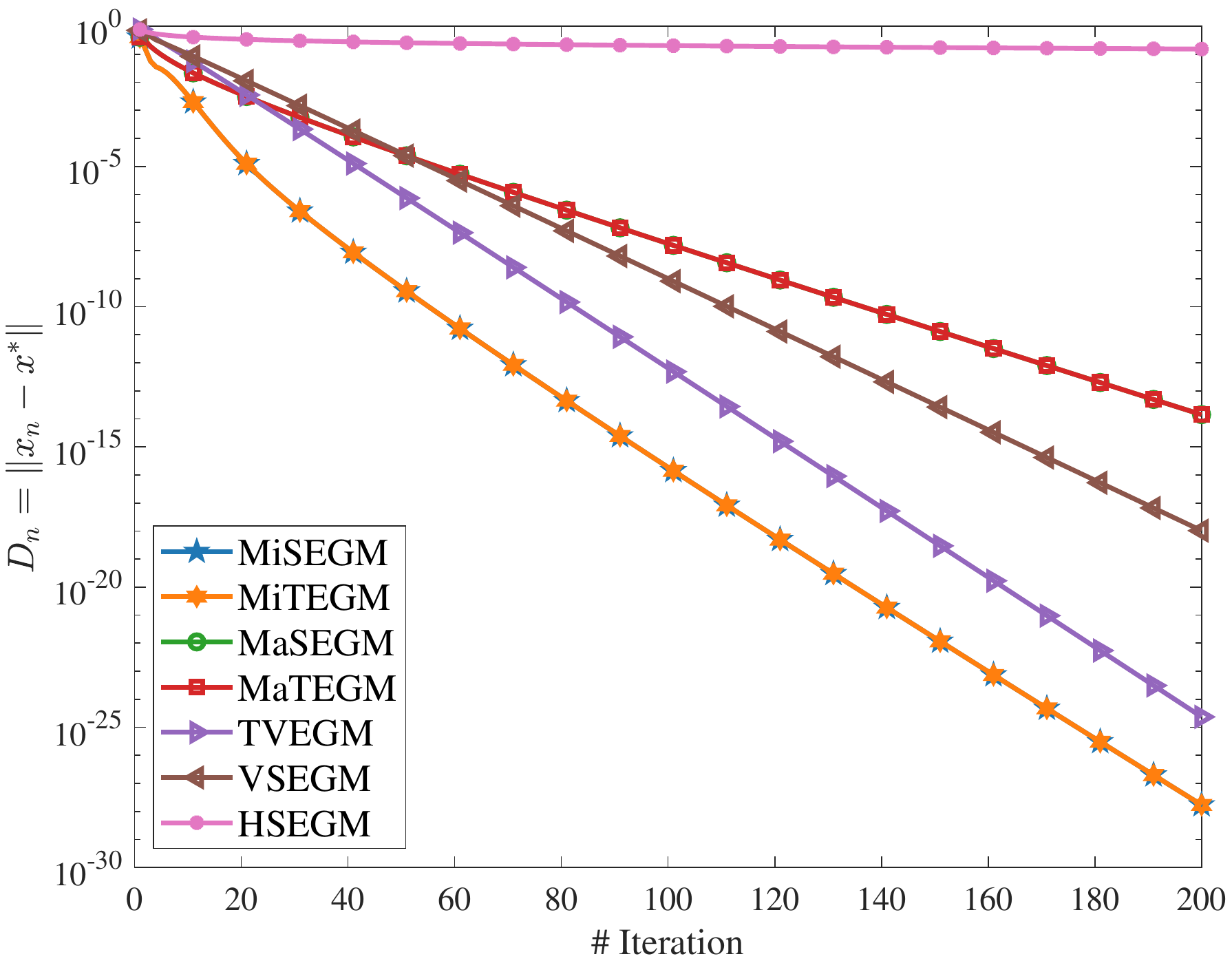}
\caption{Comparison of the number of iterations of all algorithms for Example~\ref{ex1}}
\label{ex1_fig1}
\end{figure}
\begin{figure}[htbp]
\centering
\includegraphics[scale=0.66]{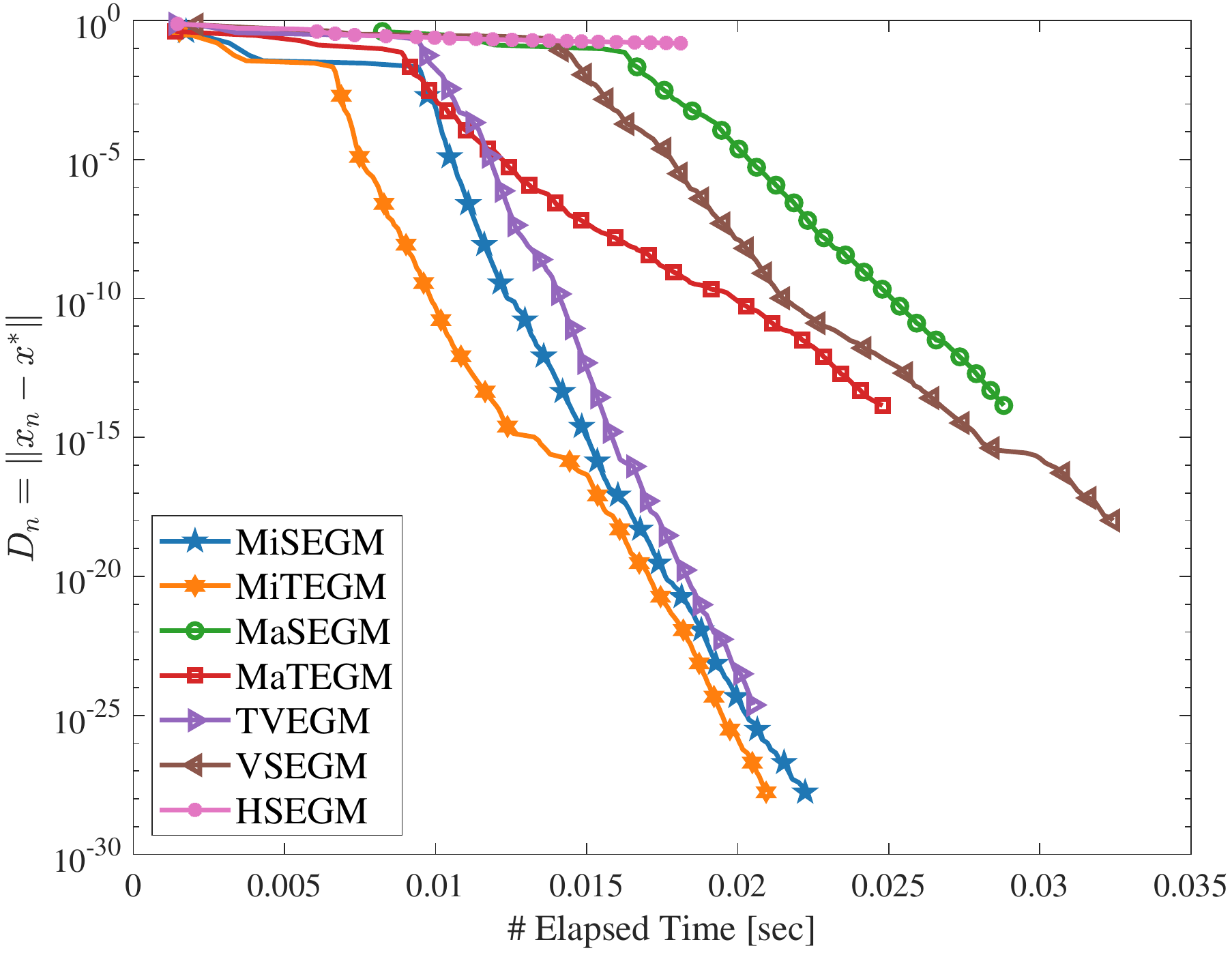}
\caption{Comparison of the elapsed time of all algorithms for Example~\ref{ex1}}
\label{ex1_fig2}
\end{figure}
\end{example}

\begin{example}\label{ex2}
Consider the linear operator $ A: R^{m}\rightarrow R^{m} $ ($ m=5 $) in the form $A(x)=Mx+q$, where $q\in R^m$ and $M=NN^{\mathsf{T}}+U+D$, $N$ is a $m\times m$ matrix, $U$ is a $m\times m$ skew-symmetric matrix, and $D$ is a $m\times m$ diagonal matrix with  its diagonal entries being nonnegative (hence $ M $ is positive symmetric definite). The feasible set $C$ is given by $C=\left\{x \in {R}^{m}:-2 \leq x_{i} \leq 5, \, i=1, \ldots, m\right\}$. It is clear that $A$ is monotone and Lipschitz continuous with constant $ L=  \|M\| $.  In this experiment, all entries of $N, D$ are generated randomly in $[0,2]$ and  $U$ is generated randomly in $[-2,2]$. Let $ q = 0 $, then the solution set is $ x^{*}=\{\mathbf{0}\} $.  The initial values $ x_{0} = x_{1} $ are randomly generated by \emph{10rand(m,1)} in MATLAB. The numerical results are shown in Figs.~\ref{ex2_fig1} and \ref{ex2_fig2}.
\begin{figure}[htbp]
\centering
\includegraphics[scale=0.66]{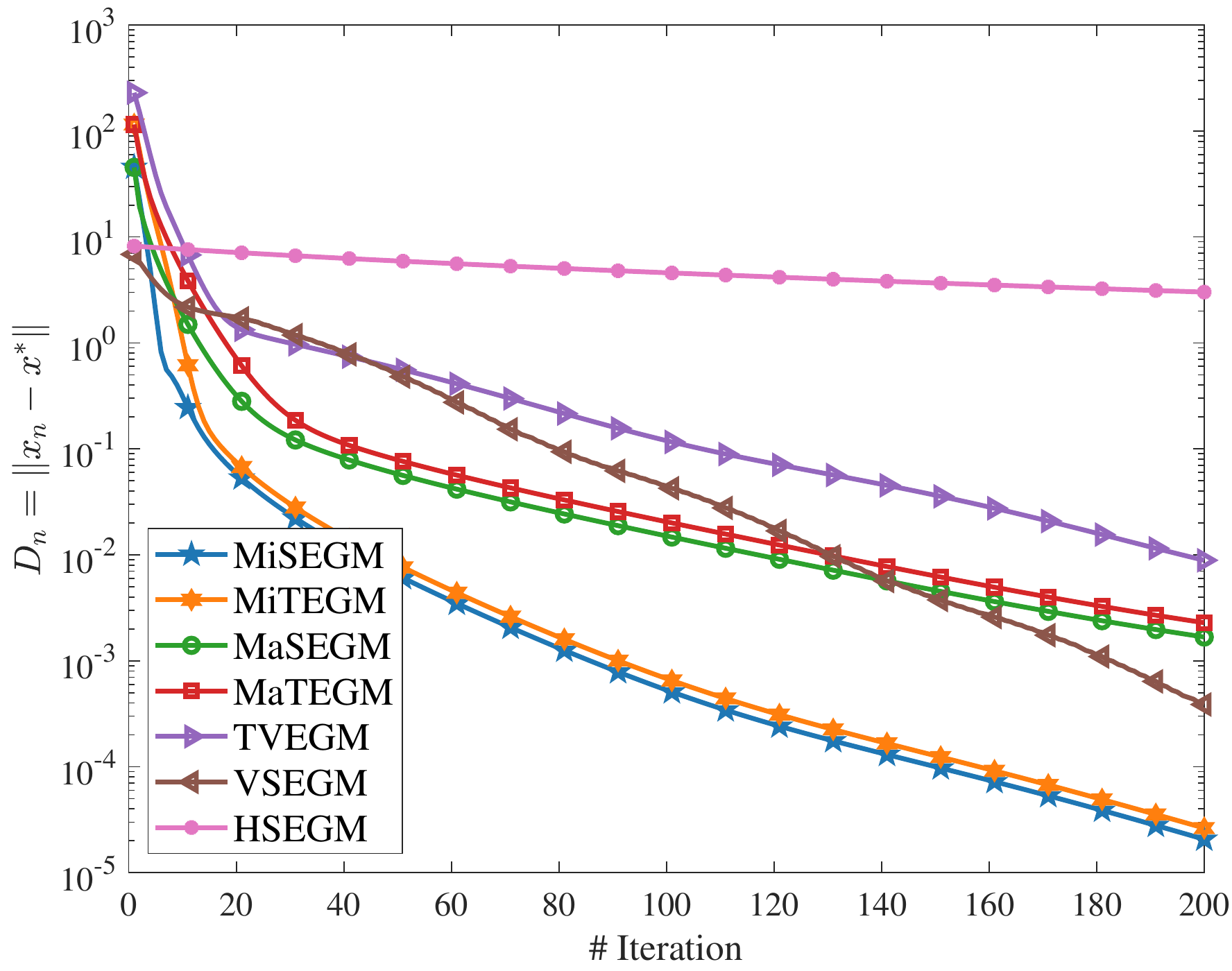}
\caption{Comparison of the number of iterations of all algorithms for Example~\ref{ex2}}
\label{ex2_fig1}
\end{figure}
\begin{figure}[htbp]
\centering
\includegraphics[scale=0.66]{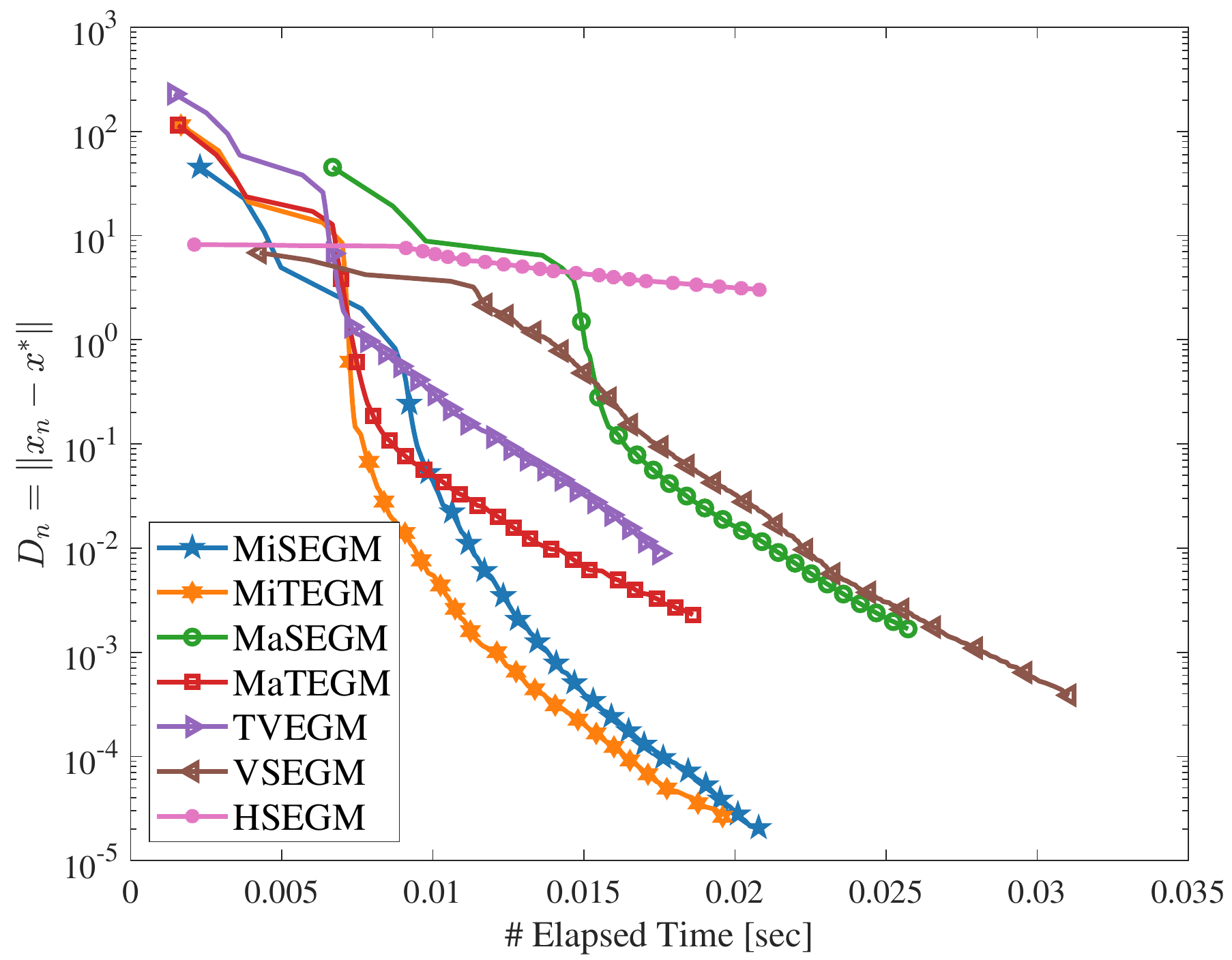}
\caption{Comparison of the elapsed time of all algorithms for Example~\ref{ex2}}
\label{ex2_fig2}
\end{figure}
\end{example}

\begin{example}\label{ex3}
Finally, we consider our problem in the Hilbert space $H=L^{2}([0,1])$ with the inner product $\langle x, y\rangle:=\int_{0}^{1} x(t) y(t) \mathrm{d} t$ and the induced norm $\|x\|:=(\int_{0}^{1}|x(t)|^{2} \mathrm{d} t)^{1 / 2},\forall x, y \in H$. Let the feasible set be the unit ball $C:=\{x \in H:\|x\| \leq 1\}$. Define an operator $A: C \rightarrow H$ by
\[
(Ax)(t)=\int_{0}^{1}\left(x(t)-G(t, s) g(x(s))\right) \mathrm{d}s +h(t), \quad  t \in[0,1],\, x \in C,
\]
where
\[
G(t, s)=\frac{2 t s \mathrm{e}^{t+s}}{e \sqrt{\mathrm{e}^{2}-1}}\,, \quad g(x)=\cos x\,, \quad h(t)=\frac{2 t \mathrm{e}^{t}}{e \sqrt{\mathrm{e}^{2}-1}}\,.
\]
It is known that $ A $ is monotone and $ L $-Lipschitz continuous with $ L = 2 $ and $ x^*(t) = \{\mathbf{0}\} $ is the solution of the corresponding variational inequality problem. Note that the projection on $ C $ is inherently explicit, that is,
\[
P_{C}(x)=\left\{\begin{array}{ll}
\frac{x}{\|x\|_{L^{2}}}, & \text { if }\|x\|_{L^{2}}>1\,; \\
x, & \text { if }\|x\|_{L^{2}} \leq 1\,.
\end{array}\right.
\]
We choose the maximum iteration of $ 50 $ as a common stopping criterion. Figs.~\ref{ex3_fig1}~and~\ref{ex3_fig2} show the numerical behaviors of all the algorithms with the starting points $ x_{0}(t)=x_{1}(t)=10e^{t} $.
\begin{figure}[htbp]
\centering
\includegraphics[scale=0.66]{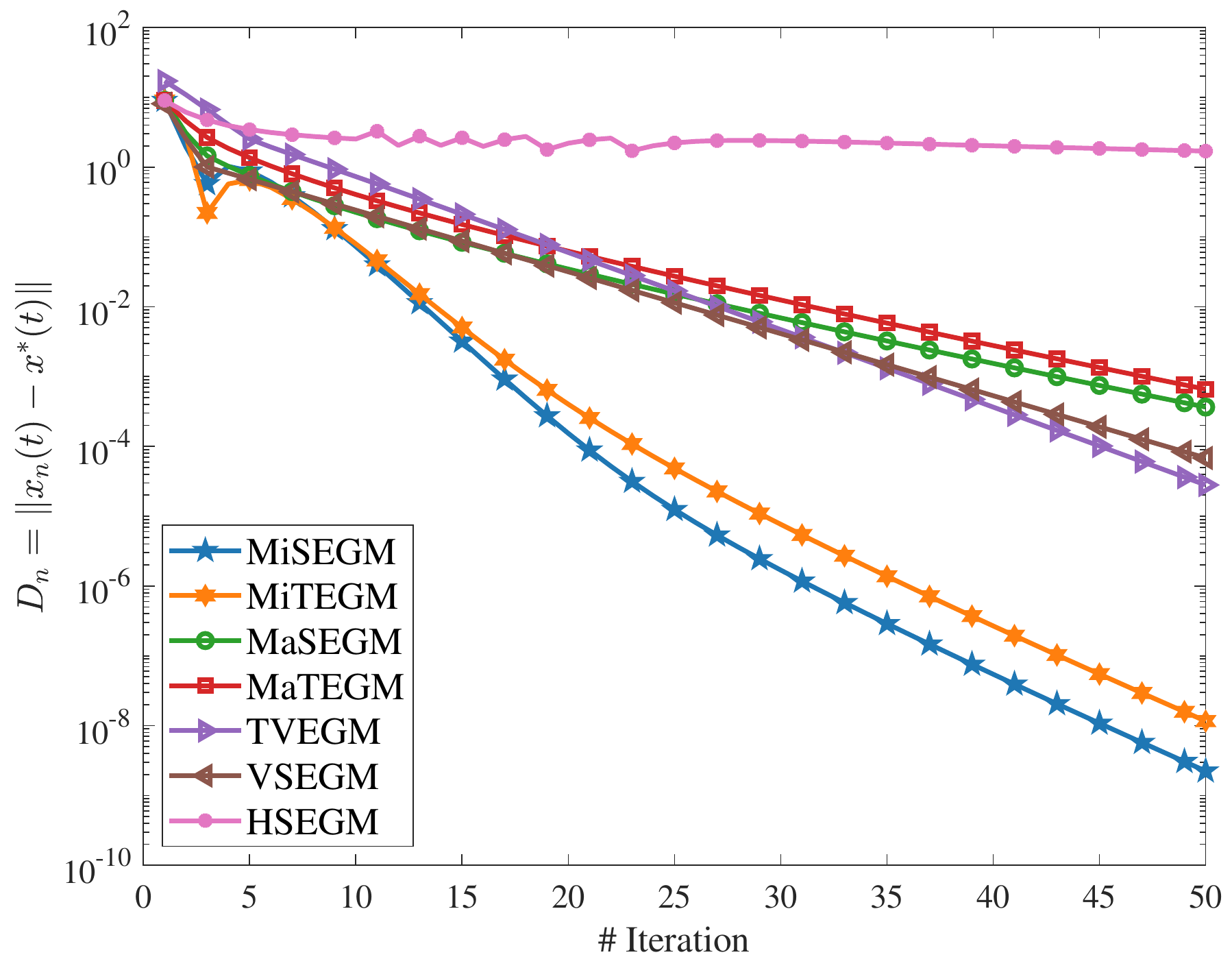}
\caption{Comparison of the number of iterations of all algorithms for Example~\ref{ex3}}
\label{ex3_fig1}
\end{figure}
\begin{figure}[htbp]
\centering
\includegraphics[scale=0.66]{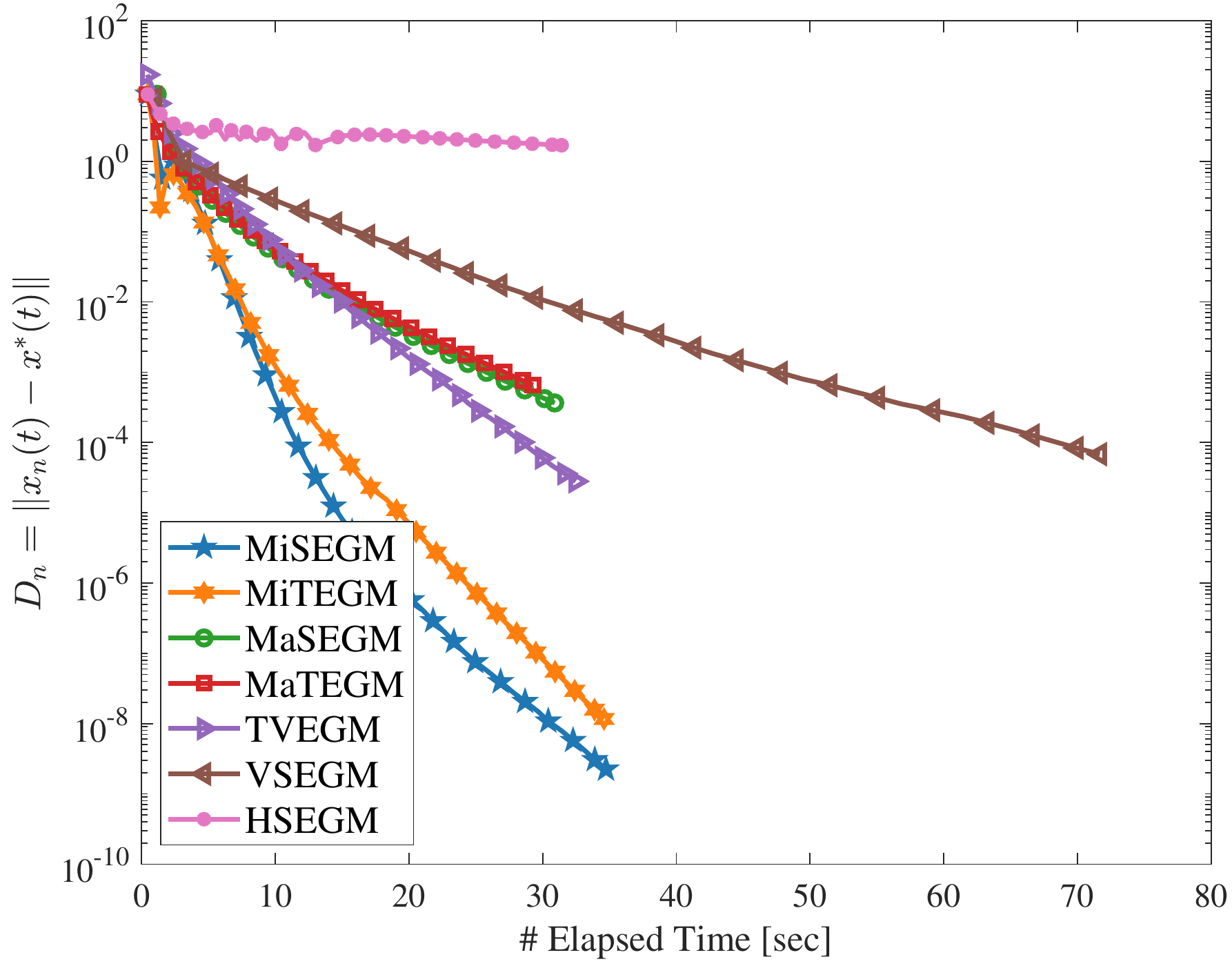}
\caption{Comparison of the elapsed time of all algorithms for Example~\ref{ex3}}
\label{ex3_fig2}
\end{figure}
\end{example}

\begin{remark}
\begin{enumerate}
\item From Figs. \ref{ex1_fig1}--\ref{ex3_fig2}, we know that our proposed algorithms outperformance the existing algorithms in terms of the number of iteration and the elapsed time.
\item It is worth noting that our algorithms converge very quickly, and there are still some oscillations because the inertial selection is too large.
\item The maximum number of iterations we choose is only $ 200 $. Note that the iteration error of algorithm~\eqref{HSEGM} is very big. In actual applications, it may require more iterations to meet the accuracy requirements.
\item We point out that since the algorithm~\eqref{VSEGM} uses the Armijo-like step size rule, which leads to taking more execution time (cf. Fig. \ref{ex3_fig2}).
\end{enumerate}
\end{remark}
\section{Conclusion}\label{sec5}
In this paper, we presented two new  iterative extragradient algorithms with a new step size for finding the solution set of a monotone, Lipschitz-continuous variational inequality problems in real Hilbert spaces. We have proved convergence theorems of the proposed algorithms under some mild conditions imposed on parameters. Some numerical examples of finite and infinite dimensions have been performed to illustrate the performance of the algorithms and compare them with previously known ones. The two algorithms obtained in this paper improve and extend the results of some existing literature.

\end{document}